\documentclass[conf]{new-aiaa}
\usepackage[utf8]{inputenc}

\usepackage{amsmath}
\usepackage{amssymb,amsfonts}
\usepackage{algorithmic}
\usepackage{graphicx}
\usepackage{textcomp}
\usepackage{xcolor}
\usepackage{amsmath}
\usepackage[version=4]{mhchem}
\usepackage{siunitx}
\usepackage{longtable,tabularx}
\usepackage{derivative}
\usepackage{physics}
\def\0{\bf \0}

\def\I{{\bf I}}

\def\0{{\bf 0}}

\def\R{\mathbb{R}}
\def\S{{\bf S}}
\def\T{{\bf T}}

\def\X{{\bf X}}

\def\Z{{\bf Z}}
\def\a{{\bf a}}
\def\b{{\bf b}}
\def\c{{\bf c}}

\def\e{{\bf e}}

\def\g{{\bf g}}
\def\h{{\bf h}}

\def\k{{\bf k}}

\def\s{{\bf s}}

\def\v{{\bf v}}
\def\w{{\bf w}}
\def\x{{\bf x}}
\def\y{{\bf y}}
\def\z{{\bf z}}

\def\T{{\rm T}}

\newtheorem{algorithm}{Algorithm}[section]

\newtheorem{theorem}{Theorem}[section]
\newtheorem{lemma}{Lemma}[section]

\title{Combining Reinforcement Learning with Arc-search Interior-Point Method for Path Planning\thanks{This paper is based on work performed at the Air Force Research Laboratory (AFRL) Control Science Center and is supported by AFOSR LRIR \#24RQCOR002 (funded by Dr. Frederick Leve). DISTRIBUTION STATEMENT A. Approved for  public release: distribution is unlimited AFRL-2026-1808 Cleared 13 APR 2026.}}

\author{Yaguang Yang\footnote{Department of Electrical and Computer Engineering Hampton University, Hampton, VA}}
\affil{Hampton University, Hampton, VA 23669, USA}

\author{Qiang Le\footnote{Department of Electrical and Computer Engineering Hampton University, Hampton, VA}}
\affil{Hampton University, Hampton, VA 23669, USA}

\author{Isaac E. Weintraub\footnote{Control Science Center, Air Force Research Laboratory, AIAA Associate Fellow. This paper is based on work performed at the Air Force Research Laboratory (AFRL) Control Science Center and is supported by AFOSR LRIR \#24RQCOR002 (funded by Dr. Frederick Leve). DISTRIBUTION STATEMENT A. Approved for  public release: distribution is unlimited AFRL-2026-1808 Cleared 13 APR 2026.}}
\affil{Air Force Research Laboratory, Wright-Patterson AFB, OH 45433, USA.}

\begin{document}

\maketitle    

\begin{abstract}
Path planning in environments containing obstacles has numerous practical applications. The problem is challenging because it is inherently nonlinear and nonconvex. Consequently, a variety of techniques have been developed to address this problem, among which machine learning and optimal control (or optimization) have emerged as two prominent approaches. In general, machine learning methods do not require a high-fidelity model, and a trained agent can often generate a feasible path in real time. However, the resulting path is not necessarily optimal with respect to performance objectives such as minimizing path length or travel time. In contrast, optimal control and optimization methods typically rely on high-fidelity models and often require computational effort that may not satisfy real-time constraints. Nevertheless, these methods are more likely to produce optimal or near-optimal solutions. To overcome the limitations of each approach while exploiting their respective strengths, this paper proposes a framework that combines reinforcement learning with an arc-search interior-point method for path planning. Numerical simulations demonstrate that the proposed approach effectively integrates the real-time decision-making capability of reinforcement learning with the optimization performance of the arc-search interior-point method, resulting in improved path-planning performance.
\end{abstract}


 
\section{Introduction}
\label{sec:introduction}
Path planning is an active area of research because of its wide range of applications~\cite{avmr25,ylw26}. In real-world scenarios, path-planning problems are often subject to restricted zones or obstacles~\cite{avmr25,lyw26a,wvchf22}. Consequently, these problems are generally not only nonlinear but also nonconvex, making them intrinsically more challenging than linear or convex optimization problems. Therefore, a variety of methods have been developed to address such problems~\cite{ylw26}.

Among the available approaches, two classes of methods are particularly attractive because of their respective advantages. The first class is based on optimal control and optimization techniques, including nonlinear optimal control~\cite{plac19}, adaptive control~\cite{wmnl24}, model predictive control~\cite{fs19}, and pseudospectral methods~\cite{lw26}. The second class is based on machine learning techniques, particularly reinforcement learning methods such as Q-learning~\cite{cx18}, Deep Q Networks (DQN)~\cite{gzlsx23}, and Deep Deterministic Policy Gradient (DDPG)~\cite{gylwc23}.

Optimal control methods typically rely on high-fidelity system models and, when successful, can produce optimal solutions. However, the computational effort required to solve the resulting optimization problems is often substantial, making these methods less suitable for real-time applications. In contrast, reinforcement learning methods generally do not require high-fidelity models. Although significant computational resources may be needed during the training phase, a trained agent can generate feasible paths much faster than optimization-based methods, making reinforcement learning particularly attractive for real-time applications~\cite{lyw26a}. Nevertheless, reinforcement learning methods typically produce feasible rather than optimal solutions. Furthermore, there may exist initial conditions for which a trained agent is unable to generate a feasible path~\cite{lyw26a}.

It is worthwhile to point out that the efficiency and effectiveness of optimal-control and optimization algorithms depend heavily on the selection of an initial point\footnote{Throughout the remainder of this paper, we distinguish between a \emph{starting point} of a path and an \emph{initial point} of a nonlinear optimization algorithm. The former refers to the initial location of a path, whereas the latter refers to the initial guess for the optimization variables. Given a starting point and an initial point, a candidate path can be constructed.}. In many cases, finding a feasible initial point may itself be a challenging task. However, if a feasible path\footnote{Throughout this paper, a path is said to be \emph{feasible} if it does not enter any restricted region and reaches the destination within a prescribed tolerance. A path is said to be \emph{optimal} if it is the shortest feasible path and reaches the destination exactly. The \emph{feasible set} is defined as the set of starting points for which a trained agent or a given optimization algorithm is able to generate a feasible path, while the \emph{optimal set} is defined as the set of starting points for which the trained agent or the algorithm is able to generate an optimal path.} generated by a machine-learning method is used as the initial point of an optimization algorithm, the efficiency and effectiveness of the optimization process can be significantly improved. Moreover, the quality of the path obtained by reinforcement learning may be enhanced because a feasible path can potentially be refined into an optimal path through subsequent optimization. More importantly, the application of an arc-search interior-point method (IPM) may enlarge both the feasible set and the optimal set. For applications with stringent real-time requirements, feasible or optimal paths corresponding to all relevant starting points may be computed offline and stored in a database. During real-time operation, the appropriate optimized path can then be retrieved directly once a starting point is specified.

Motivated by these considerations, this paper proposes a framework that combines machine learning with optimization techniques for path planning in the presence of obstacles. First, the problem is formulated as a nonlinear optimization problem tailored to the scenario under consideration. Although the present work considers only basic engagement zones (BEZ)~\cite{vw24}, the proposed framework can be readily extended to more sophisticated weapon engagement zone (WEZ) models~\cite{dzwv23}.

Three optimization algorithms are employed to solve the resulting nonlinear optimization problem: the MATLAB implementation of the interior-point method (MATLAB-IPM)\footnote{MATLAB's {\tt fmincon} implements an interior-point algorithm based on the barrier methods developed in \cite{bgn00,bmhn99}, together with line-search and trust-region techniques described in \cite{wmno06} and documented in \cite{matlabCNOA}.}, the MATLAB implementation of sequential quadratic programming (MATLAB-SQP)\footnote{MATLAB's {\tt fmincon} implements a sequential quadratic programming algorithm based on \cite{han77,power78}, incorporating a presolver \cite{cv01}, a quadratic programming solver \cite{gmsw84}, step-size strategies \cite{han77,ms83}, and additional heuristic enhancements \cite{bss88,steihaug83}, as documented in \cite{matlabCNOA}.}, and the recently developed arc-search IPM algorithm proposed in~\cite{yang25}.

Feasible paths for two path-planning problems, one involving a single obstacle and the other involving three obstacles, were previously obtained using DDPG in~\cite{lyw26a}. These paths are used as initial paths for MATLAB-IPM, MATLAB-SQP, and the arc-search IPM (ASIPM). When an initial path generated by DDPG enters a restricted region and therefore fails to satisfy the interior-point requirement, a simple heuristic is employed to construct an interior feasible path. Using the same set of initial paths, feasible-set maps and optimal-set maps are generated for all three optimization algorithms. The resulting comparisons demonstrate the advantages of combining DDPG with the arc-search IPM algorithm relative to the MATLAB-IPM and MATLAB-SQP approaches.

The remainder of the paper is organized as follows. Section~\ref{sec:formulation} presents the mathematical formulation of the path-planning problem. Section~\ref{sec:optSolution} briefly reviews the arc-search IPM algorithm. Section~\ref{sec:Testing} discusses the implementation details and numerical results. Finally, Section~\ref{sec:conclusions} summarizes the main conclusions of the paper.

\section{The problem formulation}
\label{sec:formulation}

Throughout the remainder of this paper, bold uppercase letters denote matrices, bold lowercase letters denote vectors, and ordinary letters denote scalars. For example, $\X$ denotes a matrix, $\x$ denotes a vector, and $x_{k}$ denotes the $k$-th element of $\x$. To simplify the notation, the column vector $[x_1, \ldots, x_n]^{\T}$ is written as $(x_1, \ldots, x_n)$, and the stacked vector $[\x^{\T}, \y^{\T}]^{\T}$ is written as $(\x, \y)$. Finally, given a vector $\z$, $\Z$ denotes a diagonal matrix whose diagonal entries are the components of $\z$.

Let $(\bar{x}_k,\bar{y}_k)$ denote the position of a moving object (hereafter referred to as a vehicle or an agent) at the discrete time instant $t_k$. Without loss of generality, let $(\bar{x}_0,\bar{y}_0)$ denote the initial position of the vehicle and $(\bar{x}_f,\bar{y}_f)$ denote its desired destination. Both positions are assumed to be specified prior to the path-planning process. The path from the initial position $(\bar{x}_0,\bar{y}_0)$ to the destination $(\bar{x}_f,\bar{y}_f)$ is divided into $f$ segments of equal length. Let $\theta_k$, $k=0,\ldots,f-1$, denote the vehicle heading at the beginning of the $k$-th segment, and let $r$ denote the length of each segment. The objective of the path-planning problem is to minimize the total path length while avoiding obstacles. Since the straight-line distance between the initial position and the destination is $\sqrt{(\bar{x}_f-\bar{x}_0)^2+(\bar{y}_f-\bar{y}_0)^2},$ the total path length must satisfy $f \cdot r \geq \sqrt{(\bar{x}_f-\bar{x}_0)^2+(\bar{y}_f-\bar{y}_0)^2}.$ Consequently, $r \geq \frac{\sqrt{(\bar{x}_f-\bar{x}_0)^2+(\bar{y}_f-\bar{y}_0)^2}}{f}.$

\begin{figure}[htb]
\centerline{\includegraphics[height=8cm,width=8cm]
{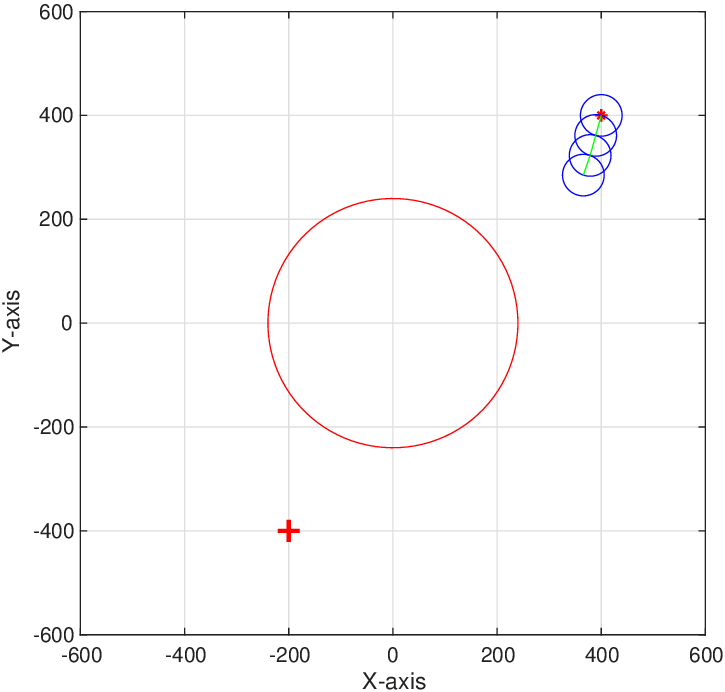}}
\caption{Path trajectory design.}
\label{pPlanning}
\end{figure}

Referring to Figure~\ref{pPlanning}, the red '*' denotes the starting point and the red '+' denotes the destination. The large red circle, centered at $(a,b)$ with radius $R$, represents the obstacle (or basic engagement zone~\cite{wvchf22}) that the vehicle is not allowed to enter\footnote{This simple obstacle model covers several practical scenarios, including certain pursuit--evasion problems described in \cite{vw24}. Furthermore, the methodology developed in this paper can be readily extended to more general constraint regions, such as those considered in \cite{dzwv23,wvchf22}.}. Starting from $(\bar{x}_0,\bar{y}_0)$, the vehicle travels a distance $r$ in the direction specified by theheading angle of $\theta_0$ and arrive at $(\bar{x}_1,\bar{y}_1)$. It then travels another distance $r$ in the direction specified by $\theta_1$ and arrive at $(\bar{x}_2,\bar{y}_2)$, and the process continues in the same manner until $(\bar{x}_f,\bar{y}_f)$ is reached. The small blue circles are centered at the starting point of every $(\bar{x}_i, \bar{y}_i)$, $i=0, 1, 2, \cdots,$ and each has radius $r$. Since the vehicle moves a fixed distance $r$ in every step, the position of the vehicle can be determined recursively from the heading angles. Given the initial position $(\bar{x}_0,\bar{y}_0)$, the following equations hold:
\begin{subequations}
\begin{gather}
\bar{x}_1 = \bar{x}_0 + r \cos(\theta_{0}), \hspace{0.1in} \bar{y}_1=\bar{y}_0 + r \sin(\theta_{0})\\
\bar{x}_k=\bar{x}_{k-1} + r \cos(\theta_{k-1}) =\bar{x}_0 +r \sum_{i=0}^{k-1} \cos(\theta_i), 
\hspace{0.1in} k=1, \ldots f-1,  \\
\bar{y}_k= \bar{y}_{k-1} + r \sin(\theta_{k-1}) =\bar{y}_0 + r \sum_{i=0}^{k-1} \sin(\theta_i), 
\hspace{0.1in} k=1, \ldots f-1, \\
\bar{x}_f=\bar{x}_0 +r \sum_{i=0}^{f-1} \cos(\theta_i), \hspace{0.1in} \bar{y}_f=\bar{y}_0 + r \sum_{i=0}^{f-1} \sin(\theta_i).
\end{gather}
\end{subequations}
Since the vehicle is not allowed to enter the obstacle region $(\bar{x}_k-a)^2+( \bar{y}_k-b)^2 \leq R^2$, and since both the initial position and the destination are assumed to lie outside the obstacle, every point along a feasible path must satisfy
\begin{equation}
\left[\bar{x}_0 + r \sum_{i=0}^{k-1} \cos(\theta_i) -a \right]^2 
+ \left[ \bar{y}_0 +r \sum_{i=0}^{k-1} \sin(\theta_i) -b \right]^2 \geq R^2, \hspace{0.1in} k=1,\ldots, f-1.
\label{cConst}
\end{equation}
To restrict the vehicle's operating region, an artificial boundary is imposed.
\begin{equation}
\left[\bar{x}_0 + r \sum_{i=0}^{k-1} \cos(\theta_i) -a \right]^2 
+ \left[ \bar{y}_0 +r \sum_{i=0}^{k-1} \sin(\theta_i) -b \right]^2 \leq D^2, \hspace{0.1in} k=1,\ldots, f-1.
\label{uConst}
\end{equation}
To avoid abrupt changes in direction, the difference in heading between two consecutive path segments is constrained by $| \theta_k -\theta_{k-1} | \leq 0.5$ radians for all $k \in \{1,\ldots,f-1\}$.

Based on the above discussion, the path-planning problem considering a single circular obstacle (BEZ) is formulated as the following optimization problem.
\begin{subequations}
\begin{alignat}{2}
\min & ~~ r \label{obj1} \\
s.t.  &   ~~ \bar{x}_f=\bar{x}_0 + r \sum_{i=0}^{f-1} \cos(\theta_i),   \label{xCoord1}  \\
&   ~~  \bar{y}_f= \bar{y}_{0} +r \sum_{i=0}^{f-1} \sin(\theta_i ),   \label{yCoord1}  \\
&  ~~   \left[\bar{x}_0 + r \sum_{i=0}^{k-1} \cos(\theta_i) -a \right]^2 
+ \left[\bar{y}_{0} + r \sum_{i=0}^{k-1} \sin(\theta_i) -b \right]^2 -R^2  \geq 0, \hspace{0.1in} k=1,\ldots, f-1,   \label{cConst1} \\
&  ~~  D^2- \left[\bar{x}_0 + r \sum_{i=0}^{k-1} \cos(\theta_i) -a \right]^2 
- \left[\bar{y}_{0} + r \sum_{i=0}^{k-1} \sin(\theta_i) -b \right]^2 \geq 0, \hspace{0.1in} k=1,\ldots, f-1,   \label{cConst2} \\
&  ~~    \theta_k - \theta_{k-1} + 0.5 \geq 0, \hspace{0.1in} k = 1,\ldots, f-1,  \label{bound1}  \\
&  ~~    \theta_{k-1} - \theta_{k} + 0.5 \geq 0, \hspace{0.1in} k = 1,\ldots, f-1,  \label{bound2} \\
&  ~~  r \geq \frac{\sqrt{(\bar{x}_f- \bar{x}_{0})^2+  (\bar{y}_f - \bar{y}_{0})^2}}{f},  \label{segLengthL} 
\end{alignat}
\label{optP1}
\end{subequations}
The problem consists of $2$ equality constraints, given by (\ref{xCoord1}) and (\ref{yCoord1}); $2f-2$ nonlinear inequality constraints, where (\ref{cConst1}) represents $f-1$ nonlinear inequality constraints, (\ref{cConst2}) represents $f-1$ nonlinear inequality constraints; $2f-2$ linear inequality constraints defined by (\ref{bound1}) and (\ref{bound2}); and $1$ boundary constraint defined by (\ref{segLengthL}). There are $n=f+1$ optimization variables $\x = (r, \theta_0, \theta_1, \ldots, \theta_{f-1})$. For convenience, the equality constraints are grouped into a vector-valued function $\h(\x)=\0$ whose components are (\ref{xCoord1}) and (\ref{yCoord1}). Similarly all inequality constraints are collected into a vector-valued inequality $\g(\x) \geq \0$ whose components are (\ref{cConst1}),  (\ref{cConst2}), (\ref{bound1}), (\ref{bound2}),  (\ref{segLengthL})). Therefore, Problem (\ref{optP1}) can be expressed in a standard form:
\begin{subequations}
\begin{alignat}{2}
\min & ~~ x_1 \label{obj2} \\
s.t.  &   ~~ \h(\x)=\0,   \label{eqC2}  \\
&  ~~   \g(\x) \geq \0.  \label{ineqC2}
\end{alignat}
\label{optP}
\end{subequations}
Denote $f(\x)=x_1 $, by introducing a relaxation vector $\s \geq \0$, Problem (\ref{optP}) can be rewritten as:
\begin{subequations}
\begin{alignat}{2}
\min & ~~ f(\x) \label{obj} \\
s.t.  &   ~~ \h(\x)=\0,   \label{eqC}  \\
&  ~~   \g(\x) + \s = \0, ~~ \s \geq \0, \label{ineqC}
\end{alignat}
\label{opt}
\end{subequations}
where $f(\x): \R^n \rightarrow \R$, $\h(\x):  \R^n \rightarrow \R^m$, and $\g(\x):  \R^n \rightarrow \R^p$. In the present problem, $n=f+1$, $m=2$, and $p=3f-2$.  To improve computational efficiency, the inequality constraints can be further partitioned into two groups: nonlinear inequality constraints and linear inequality constraints. Since the second-order derivatives associated with the linear inequality constraints are identically zero, their evaluation can be omitted during the computation of second-order information (see (\ref{2ndL}) and (\ref{secondOrderM})). This separation reduces the computational burden without affecting the optimization results.

\section{Optimal path planning using arc-search IPM}
\label{sec:optSolution}

We adopt an arc-search interior-point method (IPM) proposed in \cite{yang25} to solve the nonlinear optimization problem (\ref{opt}). The Lagrangian of Problem (\ref{opt}) is defined as
\begin{equation}
L(\x,\y,\w,\s,\z)=f(\x)-\y^{\T}\h(\x)-\w^{\T}(\g(\x)-\s)-\z^{\T}\s,
\label{lagrangian1}
\end{equation}
where $\y \in \R^m, \w>\0 \in \R^p_+$, and $\z>0 \in \R^p_+$  are multipliers. Let $\v= (\x,\y,\w,\s,\z) \in \R^{n+m+3p}$ denote a stacked column vector composed of primal and slack variables, and multipliers. The gradient of the Lagrangian function with respect to $\x$ and $\s$ are given by
\begin{equation}
\nabla_{\x} L(\v)=\nabla f(\x)-\nabla \h(\x) \y-\nabla \g(\x) \w,
\label{dlagrangian}
\end{equation}
and the Hessian matrix of $L$ is given by
\begin{equation}
\nabla_{\x}^2 L(\v) = \nabla_{\x}^2 f(\x) 
- \sum_{i=1}^m \nabla_{\x}^2 h_i(\x)  y_i
- \sum_{i=1}^p \nabla_{\x}^2 g_i(\x)  w_i.
\label{2ndL}
\end{equation}

Most optimization algorithms find the optimal solution by searching a KKT point \cite{nw06}, which is equivalent to solve a system of nonlinear equations:
\begin{align}
	\k(\v)  = \left[ \begin{array}{c}
		\nabla_{\x} L(\v)  \\
		\h(\x)  \\
		\g(\x)-\s \\
		\w-\z  \\
		\Z\s
	\end{array} \right] =  \0,  \
	(\w, \s, \z) \in \R_+^{3p}, 
	\label{KKT1}
\end{align}
where $\k:\R^{n+m+3p} \to \R^{n+m+3p}$. A sufficient condition for KKT to hold is 
\begin{equation}
\phi(\v)=\k(\v)^{\T} \k(\v) =0.
\end{equation}

Taking the first-order derivative of (\ref{KKT1}) yields
\begin{equation}
\left[ \begin{array}{ccccc}
\nabla_{\x}^2 L({\v^k}) & -\nabla \h({\x^k}) & -\nabla \g({\x^k}) & \0  & \0 \\
\left( \nabla \h({\x^k})\right)^{\T} & \0 & \0 & \0 & \0  \\
\left( \nabla \g({\x^k})\right)^{\T} & \0 & \0 & -\I  & \0 \\
\0  & \0  &  \I  &  \0  &  -\I  \\
\0 & \0  & \0 & \Z^k & \S^k
\end{array} \right]
\left[ \begin{array}{c}
\dot{\x} \\ \dot{\y} \\ \dot{\w} \\ \dot{\s}  \\ \dot{\z}
\end{array} \right]
= \left[ \begin{array}{c}
\nabla_{\x} L({\v^k}) \\
\h({\x^k}) \\
\g({\x^k})-{\s^k} \\
{\w^k}-{\z^k}  \\
\Z^k {\s^k}  -\sigma {\mu}_k \e
\end{array} \right].
\label{firstOrderM}
\end{equation}
where $\e$ is an all-one vector of dimension $p$. Taking the derivative of (\ref{firstOrderM}) and ignoring the third-order derivative yields
\footnotesize
\begin{eqnarray}
\left[ \begin{array}{ccccc}
\nabla_{\x}^2 L(\v^k) & -\nabla \h(\x^k) & -\nabla \g(\x^k) & \0  & \0 \\
\left( \nabla \h(\x^k)\right)^{\T} & \0 & \0 & \0 & \0  \\
\left( \nabla \g(\x^k)\right)^{\T} & \0 & \0 & -\I  & \0 \\
\0  & \0  &  \I  &  \0  &  -\I  \\
\0 & \0 & \0  & \Z^k & \S^k
\end{array} \right]
\left[ \begin{array}{c}
\ddot{\x} \\ \ddot{\y} \\ \ddot{\w} \\ \ddot{\s} \\ \ddot{\z}
\end{array} \right]
=  \left[ \begin{array}{c}
2[(\nabla_{\x}^2 \g(\x^k))\dot{\w} \dot{\x}+(\nabla_{\x}^2 \h(\x^k))\dot{\y} \dot{\x}] \\
-(\nabla_{\x}^2 \h(\x^k))^{\T} \dot{\x} \dot{\x} \\
-(\nabla_{\x}^2 \g(\x^k))^{\T} \dot{\x} \dot{\x} \\
\0  \\
-2\dot{\Z} \dot{\s} 
\end{array} \right],
\label{secondOrderM}
\end{eqnarray}
\normalsize
where the computational formulas for the elements on the right-hand side of (\ref{secondOrderM}) are given in \cite[Appendix 1]{yiy18}:
\begin{eqnarray}
\nabla_{\x}^2 \h(\x)\dot{\y} \dot{\x} 
&=&
\frac{\partial \left( \frac{\partial  \h(\x)}{\partial \x } \dot{\y} \right) }
{\partial \x} \dot{\x} 
=  \sum_{i=1}^m \dot{y}_i \frac{\partial}{\partial \x}
\left[ \begin{array}{c}
	\frac{\partial  h_i(\x)}{\partial x_1 } \\
	\vdots \\
	\frac{\partial  h_i(\x)}{\partial x_n} 
\end{array} \right] \dot{\x}
=  \sum_{i=1}^m \dot{y}_i
\left( \nabla_{\x}^2 h_i(\x) \right)  \dot{\x} 
\label{hyx}
\end{eqnarray}
\begin{eqnarray}
	\nabla_{\x}^2 \h(\x)^{\T} \dot{\x}  \dot{\x}
	&=&
	\left( \frac{\partial \left( \left( \frac{\partial  \h(\x)}
		{\partial \x } \right)^{\T} \dot{\x} \right) }{\partial \x} \right)^{\T} \dot{\x}
	=  
	\left[ \begin{array}{c}
		\dot{\x}^{\T} \left( \nabla_{\x}^2 h_1(\x) \right)  \dot{\x} \\
		\vdots \\
		\dot{\x}^{\T} \left( \nabla_{\x}^2 h_m(\x) \right)  \dot{\x} 
	\end{array} \right] 
	\label{hxx}
\end{eqnarray}
\begin{eqnarray}
\nabla_{\x}^2 \g(\x)\dot{\w} \dot{\x} 
&=&
\frac{\partial \left( \frac{\partial  \g(\x)}{\partial \x } \dot{\w} \right) }
{\partial \x} \dot{\x} 
=  \sum_{i=1}^n \dot{w}_i \frac{\partial}{\partial \x}
\left[ \begin{array}{c}
	\frac{\partial  g_i(\x)}{\partial x_1 } \\
	\vdots \\
	\frac{\partial g_i(\x)}{\partial x_n} 
\end{array} \right] \dot{\x}
= \sum_{i=1}^n \dot{w}_i
\left( \nabla_{\x}^2 g_i(\x) \right)  \dot{\x}  
\label{gzx}
\end{eqnarray}
\begin{eqnarray}
\nabla_{\x}^2 \g(\x)^{\T} \dot{\x}  \dot{\x}
&=&
\left( \frac{\partial \left( \left( \frac{\partial \g(\x)}
	{\partial \x } \right)^{\T} \dot{\x} \right) }{\partial \x} \right)^{\T} \dot{\x}
=  
\left[ \begin{array}{c}
	\dot{\x}^{\T} \left( \nabla_{\x}^2 g_1(\x) \right)  \dot{\x} \\
	\vdots \\
	\dot{\x}^{\T} \left( \nabla_{\x}^2 g_p(\x) \right)  \dot{\x} 
\end{array} \right].
\label{gxx}
\end{eqnarray}

An observation from the fourth component of (\ref{KKT1}), (\ref{firstOrderM}), and (\ref{secondOrderM}) leads to the following lemma.

\begin{lemma}\upshape{\cite{yiy18}}
If $\w^0=\z^0$, then $\w^k=\z^k$.
\label{wEqz}
\end{lemma}

This observation is important because it reduces computational cost: it is sufficient to compute either $\w$ or $\z$ but not both. Another key idea in interior-point methods is to trace the central path \cite{megiddo89}. However, computing the central path explicitly is generally computationally impractical. To address this difficulty, arc-search techniques have been proposed, which approximate the central path by a segment of an ellipse \cite{yang20}.
\begin{equation}
	{\cal E}=\lbrace \v (\alpha): 
	\v (\alpha) =
	\vec{\a}\cos(\alpha)+\vec{\b}\sin(\alpha)+\vec{\c}, \alpha \in [0, 2\pi] \rbrace,
	\label{ellipse}
\end{equation}

\begin{theorem}\label{theorem:ellip}
	\upshape{\cite{yang20}}
	Suppose that the ellipse ${\cal E}$ defined in \textrm{(\ref{ellipse})} passes through 
	the current iterate ${\v^k}$ at $\alpha=0$, and its first- and second-order derivatives at 
	$\alpha=0$ are given by $\dot{\v}$ and $\ddot{\v}$, respectively. Then the trajectory on ${\cal E}$ can be expressed
	as ${\v}(\alpha) = ({\x}(\alpha), {\y}(\alpha), {\w}(\alpha), {\s}(\alpha), {\z}(\alpha))$, 
	which can be calculated by
	\begin{align}
		{\v}(\alpha) = {\v^k} - \dot{\v}\sin(\alpha)+\ddot{\v}(1-\cos(\alpha)).
		\label{vAlpha}
	\end{align}
	\label{arcFormula}
\end{theorem}

To determine the largest step size $\alpha_{w_i}$ that preserves the nonnegativity condition $w_i \geq 0$, a set of formulas of $\alpha_{w_i}$ was derived in \cite{yang20}. This yields the maximum admissible step size $\tilde{\alpha}_k$ such as $\w \geq \0$, given by
\begin{equation}
\tilde{\alpha}_k=\min_{i \in \lbrace 1,\ldots, p \rbrace}
\left\{ \min \left\{\alpha_{w_i}^k, \frac{\pi}{2} \right\} \right\}.
\label{alpha}
\end{equation}
Next a search is performed to determine $\bar{\alpha} \in (0, \tilde{\alpha}_k]$ such that the inequality constraints $\g(\x) > 0$ are satisfied. To ensure sufficient decrease of the merit function $\phi(\v(\alpha))$, an additional search is carried out to find $\check{\alpha} \in (0, \bar{\alpha}_k]$ such that for some $\delta_2 \geq 0$, the following inequality holds:
\begin{equation}
	\check{\alpha}_k = \max\left\{ \alpha \in \left(0, \bar{\alpha}_k \right] : 
	\phi({\v^k}(\alpha)) < \phi({\v^k}) - \delta_2 \right\}.
	\label{alphaCheck}
\end{equation}
Finally, to ensure convergence of the arc-search IPM algorithm, $\hat{m}_k(\alpha)$ is defined in \cite{yang25} as 
\begin{equation}
\hat{m}_k(\alpha)  = \min (\Z^k(\alpha) \s^k(\alpha) )- \frac{1}{2} 
\frac{\Z^0\s^0}{\phi (\v^0)} \min (\phi (\v(\alpha) )),
\label{measurePos}
\end{equation}
where the minimum is taken component-wise. A subsequent search is then performed to determineand $\hat{\alpha} \in (0, \bar{\alpha}_k]$ such that
\begin{equation}
\hat{\alpha}_k = \max\left\{ \alpha \in \left(0, \check{\alpha}_k \right] : \hat{m}_k(\alpha) \ge 0 \right\}.
\label{hatAlpha}
\end{equation}

Therefore, the arc-search IPM algorithm can be stated as follows:

\begin{algorithm}{\bf \cite{yang25}} 
\begin{algorithmic}[1]
\STATE Parameters: $\epsilon>0$, $\delta_1>0$, $\delta_2>0$,  
$\rho \in (0, \frac{1}{2})$, and $\bar{\sigma} \in [0, \frac{1}{2})$.
\STATE Initial point: $\v^0 = (\x^0,\y^0,\w^0,\s^0,\z^0)$ 
such that $(\w^0,\s^0,\z^0) > \0$, $\g(\x^0) >0$,  and $\w^0=\z^0$.
\FOR{  k=0,1,2,\ldots}
	\STATE  Calculate $\h(\x^k)$, $\g(\x^k)$, $\nabla_{\x} \h(\x^k)$, $\nabla_{\x} \g(\x^k)$, and $\nabla_{\x} L(\v^k)$. If $\phi(\v^k) \le \epsilon$, stop.
	\STATE Calculate 
     $\nabla_{\x}^2 \h(\x^k)$, $\nabla_{\x}^2 \g(\x^k)$, and $\nabla_{\x}^2 L(\v^k)$.
	\STATE Select $\sigma_k$ such that 
	$\bar{\sigma} \le \sigma_k <  \min \{ \frac{1}{8}, \phi(\v^k) p/\mu^2 \}$.
	\STATE Calculate $\dot{\v}^k$ by solving (\ref{firstOrderM}) at ${\v} = \v^k$.
	\STATE Calculate the righthand side of (\ref{secondOrderM}) using (\ref{hyx}), (\ref{hxx}), (\ref{gzx}), and (\ref{gxx}).
	\STATE Calculate $\ddot{\v}^k$ by solving (\ref{secondOrderM})  at ${\v} = \v^k$.
	\STATE Calculate $\tilde{\alpha}_k$ using (\ref{alpha}) and search $\bar{\alpha}_k \in (0, \tilde{\alpha}_k]$
     such that $\g(\x(\bar{\alpha}_k)) > \0$.
	\STATE Search $\check{\alpha}_k \in (0, \bar{\alpha}_k]$ such that (\ref{alphaCheck}) holds.
     \STATE  Determine $\hat{\alpha}_k > 0$ using (\ref{measurePos}) and    
     (\ref{hatAlpha}).
     \STATE  Set $\alpha_k = \min \{ \hat{\alpha}_k, \check{\alpha}_k \}$.
     \STATE  Update $(\x^{k+1},\w^{k+1}) =  (\x^k,\w^k)
	- (\dot{\x}^k,\dot{\w}^k) \sin(\alpha_k) 
     + (\ddot{\x}^k,\ddot{\w}^k) (1-\cos(\alpha_k)).$
     \STATE  
		Set $\s^{k+1} = \g(\x^{k+1}) >\0$ and $\z^{k+1}=\w^{k+1}$.
	\STATE $k+1 \rightarrow k$.
\ENDFOR 
\end{algorithmic}
\label{mainAlgo}
\end{algorithm}
The convergence of this algorithm is proved in \cite{yang25}.

\section{Implementation and testing results}\label{sec:Testing}

First, we present implementation details of Algorithm \ref{mainAlgo}.

\subsection{Assembling systems of equations (\ref{firstOrderM}) and (\ref{secondOrderM})}
To implement the algorithm, we need to compute $\nabla_{\x} h_k(\x)$ and $\nabla_{\x}^2 h_k(\x)$ for $k=1,2$, as well as $\nabla_{\x} g_k(\x)$ and $\nabla_{\x}^2 g_k(\x)$ for $k=1,2, \ldots, f-1$ in order to assemble the systems of equations (\ref{firstOrderM}) and (\ref{secondOrderM}). Taking the first- and second-order derivatives of (\ref{xCoord1}) and (\ref{yCoord1}) yields the following expressions:
\begin{equation}
	\nabla_{\x} h_1(\x) = \left(\sum_{i=0}^{f-1}\cos(\theta_i), -r\sin(\theta_0),-r\sin(\theta_1), \ldots, -r\sin(\theta_{f-1}) \right),
	\label{dh1}
\end{equation}
\begin{equation}
	\nabla_{\x}^2 h_1(\x) = \left[ \begin{array}{cccccc} 
		0 & -\sin(\theta_0) & -\sin(\theta_1) & \ldots  & \ldots & -\sin(\theta_{f-1}) \\ 
		 -\sin(\theta_0) & -r\cos(\theta_0) & 0 & \ldots  & \ldots & 0 \\
		 -\sin(\theta_1) & 0 & -r\cos(\theta_1) & 0 & \ldots & 0 \\
		 \vdots &  \vdots & \ddots  & \ddots  & \ddots  & 0 \\
		 -\sin(\theta_{f-2}) & 0  & \ldots &   0 & -r\cos(\theta_{f-2})  & 0 \\
		 -\sin(\theta_{f-1}) & 0 &  \ldots &  \ldots & 0 & -r\cos(\theta_{f-1}) 
		 \end{array} \right],
	\label{ddh1}
\end{equation}
\begin{equation}
	\nabla_{\x} h_2(\x) = \left(\sum_{i=0}^{f-1}\sin(\theta_i), 
	r\cos(\theta_0), r\cos(\theta_1), \ldots, r\cos(\theta_{f-1}) \right),
	\label{dh2}
\end{equation}
\begin{equation}
	\nabla_{\x}^2 h_2(\x) = \left[ \begin{array}{cccccc} 
		0 & \cos(\theta_0) & \cos(\theta_1) & \ldots  & \ldots & \cos(\theta_{f-1}) \\ 
		\cos(\theta_0) & -r\sin(\theta_0) & 0 & \ldots  & \ldots & 0 \\
		\cos(\theta_1) & 0 & -r\sin(\theta_1) & 0 & \ldots & 0 \\
		\vdots &  \vdots & \ddots  & \ddots  & \ddots  & 0 \\
		\cos(\theta_{f-2}) & 0  & \ldots &   0 & -r\sin(\theta_{f-2})  & 0 \\
		\cos(\theta_{f-1}) & 0 &  \ldots &  \ldots & 0 & -r\sin(\theta_{f-1}) 
	\end{array} \right],
	\label{ddh2}
\end{equation}
Taking the first-order derivatives of (\ref{cConst1}), for $k=1,\ldots,f-1$, yields
\begin{equation}
	\nabla_{\x} g_k(\x) = \left[ \begin{array}{c}
2	\left(r \sum_{i=0}^{k-1}\cos(\theta_i)-a \right) \left(\sum_{i=0}^{k-1}\cos(\theta_i) \right) + 2 	\left(r \sum_{i=0}^{k-1}\sin(\theta_i)-b \right) \left(\sum_{i=0}^{k-1}\sin(\theta_i) \right) \\
-2r	\left(r \sum_{i=0}^{k-1}\cos(\theta_i)-a \right) \sin(\theta_0) + 2r	\left(r \sum_{i=0}^{k-1}\sin(\theta_i)-b \right) \cos(\theta_0) \\
\vdots \\
-2r	\left(r \sum_{i=0}^{k-1}\cos(\theta_i)-a \right) \sin(\theta_j) + 2r	\left(r \sum_{i=0}^{k-1}\sin(\theta_i)-b \right) \cos(\theta_j) \\
\vdots \\
-2r	\left(r \sum_{i=0}^{k-1}\cos(\theta_i)-a \right) \sin(\theta_{k-1}) + 2r	\left(r \sum_{i=0}^{k-1}\sin(\theta_i)-b \right) \cos(\theta_{k-1}) \\
0 \\
\vdots \\
0
	\end{array} \right], 
	\label{dgk}
\end{equation}
To obtain the second-order derivatives of (\ref{cConst1}), we differentiate (\ref{dgk}) once more. This allows the entries of the Hessian matrix to be written as follows:
\begin{equation}
	\frac{\partial^2 g_k(\x)}{\partial r^2}  = 
		2 \left( \sum_{i=0}^{k-1}\cos(\theta_i) \right)^2 +2 \left( \sum_{i=0}^{k-1}\sin(\theta_i) \right)^2, 
	\label{ddgkdr2}
\end{equation}
\begin{equation}
	\frac{\partial^2 g_k(\x)}{\partial r \partial \theta_i}  = 
	-2 \left[ \left( 2r \sum_{i=0}^{k-1}\cos(\theta_i) -a \right)\sin(\theta_i) -  \left( 2r  \sum_{i=0}^{k-1}\sin(\theta_i) -b \right)\cos(\theta_i) \right], \hspace{0.1in} i=0,1,\ldots,k-1,
	\label{ddgkdr3}
\end{equation}
\begin{equation}
	\frac{\partial^2 g_k(\x)}{\partial \theta_j^2}  = 
	2r \left[ r+ \left( r \sum_{i=0}^{k-1}\cos(\theta_i) -a \right)\cos(\theta_j) - \left( r  \sum_{i=0}^{k-1}\sin(\theta_i) -b \right)\sin(\theta_j) \right], \hspace{0.1in} j=0,1,\ldots,k-1.
	\label{ddgkdjdj}
\end{equation}
\begin{equation}
	\frac{\partial^2 g_k(\x)}{\partial \theta_i \partial \theta_j}  = 
	2r^2 \sin(\theta_i) \sin(\theta_j) + 2 r^2  \cos(\theta_i) \cos(\theta_j), \hspace{0.1in}  i,j=0,1,\ldots,k-1, \hspace{0.1in} i \neq j.
	\label{ddgkdidj}
\end{equation}
The Hessian matrix is therefore given by
\begin{eqnarray}
\frac{\partial^2 g_k(\x)}{\partial \x^2} =	\left[ \begin{array}{c|r}
\begin{matrix} 
		\frac{\partial^2 g_k(\x)}{\partial r^2} & \frac{\partial^2 g_k(\x)}{\partial r \partial \theta_0} & \ldots  & \frac{\partial^2 g_k(\x)}{\partial r \partial \theta_{k-1}} \\
		\frac{\partial^2 g_k(\x)}{\partial r \partial \theta_0} & \frac{\partial^2 g_k(\x)}{ \partial \theta_0^2} & \ldots & \frac{\partial^2 g_k(\x)}{\partial \theta_{0} \partial \theta_{k-1}} \\
		\vdots & \vdots & \ddots  & \vdots \\
		\frac{\partial^2 g_k(\x)}{\partial r \partial \theta_{k-1}}  & \ldots  & \ldots  &  \frac{\partial^2 g_k(\x)}{ \partial \theta_{k-1}^2}  
\end{matrix} & \mbox{ $\0$} \\ \hline \mbox{ $\0$} &  \mbox{ $\0$}
	\end{array} \right], \hspace{0.1in} k=1, \ldots, f-1.
	\label{HessionM}
\end{eqnarray}

These formulas enable the assembly of the two linear systems of equations (\ref{firstOrderM}) and (\ref{secondOrderM}). Solving these systems yields the vectors $\dot{\v}$ an $\ddot{\v}$. With these quantities, the arc-search IPM can be implemented using Theorem \ref{arcFormula}, while ensuring the conditions $\w^k> \0$, $\s^k > \0$, $\z^k > \0$, $\g(\x^k) > \0$, together with (\ref{alphaCheck}), and (\ref{hatAlpha}), in order to determine the apppropriate step size and updated $\v^{k+1}$. This process is repeated until an optimal solution $\x^*=(r^*, \theta_0^*, \ldots, \theta_{f-1}^*)$ is obtained. The optimal path can then be constructed from the starting location and the computed optimal solution. For large-scale problems, explicitly using the formulas to calculate $\nabla_{\x} h_k(\x)$ and $\nabla_{\x}^2 h_k(\x)$ for $k=1,2$, and $\nabla_{\x} g_k(\x)$ and $\nabla_{\x}^2 g_k(\x)$ for $k=1,2, \ldots, f-1$ can be tedious, and manually coding the large matrices is prone to errors. Therefore, automatic differentiation \cite{nw06} is also employed to compute these derivatives directly from $\h(\x)$ and $\g(\x)$ for the path planning problem. On the other hand, since $\nabla_{\x}^2 h_k(\x)$ and $\nabla_{\x}^2 g_k(\x)$ for small $k$ are sparse, using explicit formulas given in (\ref{ddh1}), (\ref{ddh2}), and (\ref{HessionM}) instead of automatic differentiation can further reduce computational cost.

\subsection{Path planning in the presence of single obstacle}

\begin{figure}[htb]
\centerline{\includegraphics[height=6cm,width=7cm]{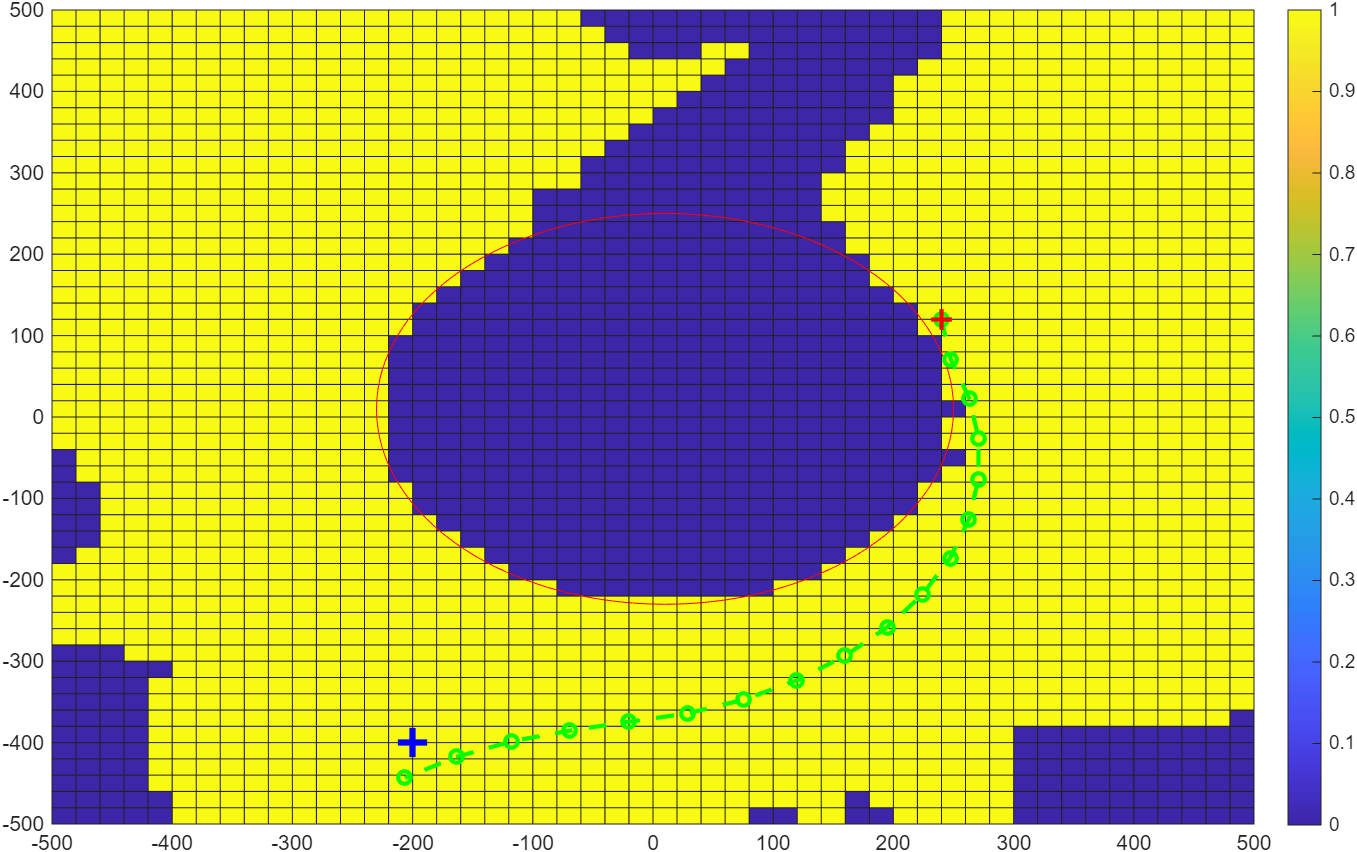}}
\caption{Feasible set found by the trained DDPG.}
\label{DDPGfeasible}
\end{figure}

\begin{figure}[htb]
\centerline{\includegraphics[height=5cm,width=8cm]{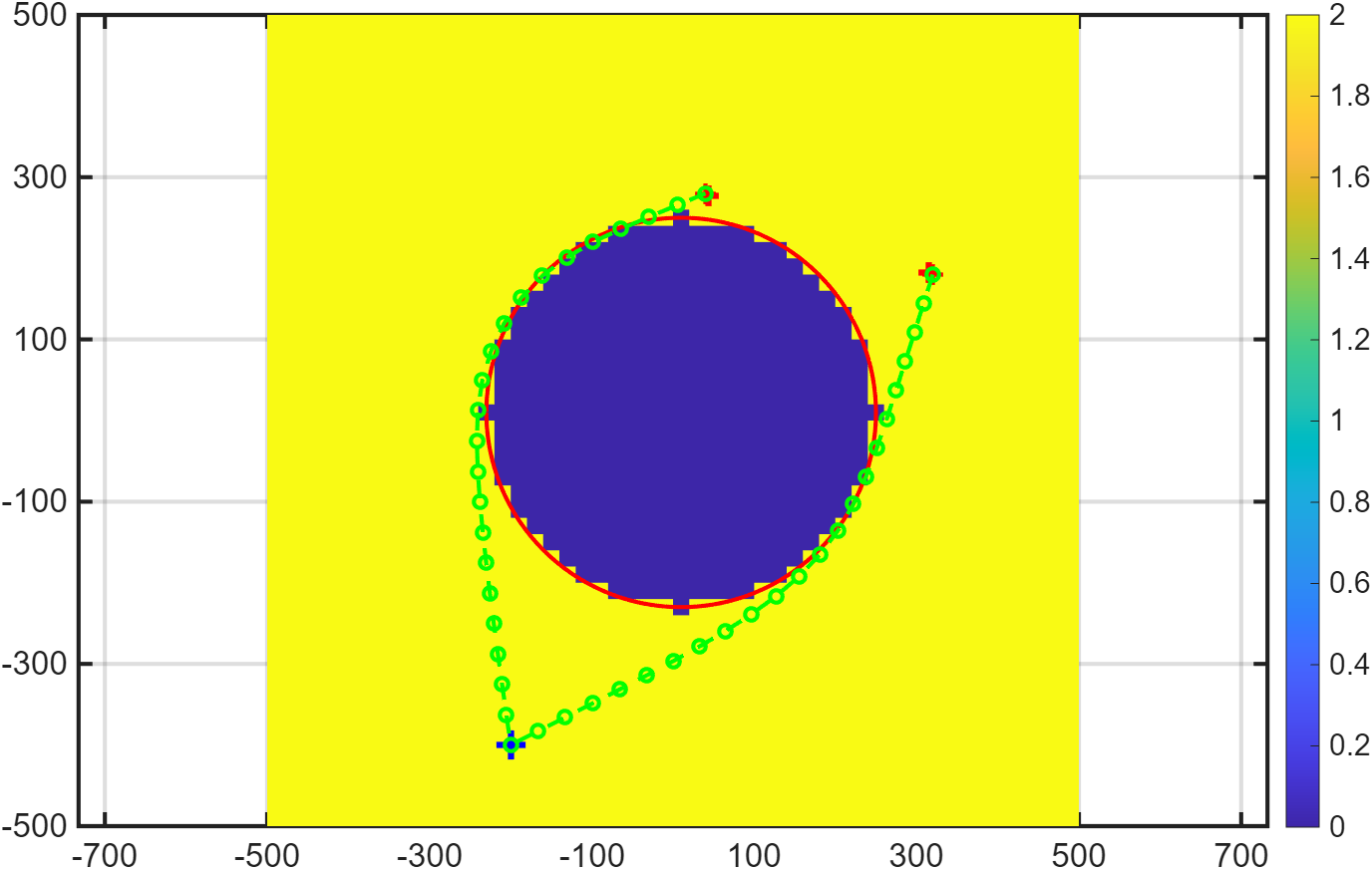}}
\caption{Optimal set found by DDPG-arc-search IPM.}
\label{arcSearch_feasible}
\end{figure}

We have shown in the previous sections that the path-planning problem can be formulated as the optimization problem (\ref{optP}). We propose using an arc-search interior-point method (IPM) to solve this problem. It is well known that a good initial path is important for optimization algorithms to efficiently find an optimal solution; however, there is no universal approach for generating such a path.

On the other hand, feasible paths can be generated using machine-learning methods, such as Deep Deterministic Policy Gradient (DDPG). In \cite{lyw26a}, DDPG is employed to identify a set of starting points from which the agent can safely approach the destination without entering any obstacle region. We refer to this set of starting points as the {\it feasible set}. A notable advantage of DDPG is that, once trained, the agent can generate a feasible path from any starting point in the feasible set to the destination in real time.

However, directly using a DDPG agent has two limitations. First, the resulting feasible path is generally not optimal. Second, the trained agent may fail to find a feasible path for starting points outside the feasible set. In this section, we combine machine-learning techniques with the arc-search IPM to exploit the strengths of both approaches. Specifically, the machine-learning agent is used to generate feasible paths, while the arc-search IPM is employed to efficiently compute optimal paths offline. For starting points for which the DDPG agent cannot generate a feasible path, the arc-search IPM is used to construct an initial path and search for an optimal solution. If an optimal solution cannot be obtained within the prescribed computational limits but a feasible solution is found, the feasible solution is also stored for subsequent real-time decision making. Consequently, the computational burden of optimization is shifted to the offline stage and does not affect real-time applications.

We define the {\it optimal set} as the set of starting points for which an optimization algorithm successfully computes optimal paths. More generally, we define the {\it feasible set} as the set of starting points for which either the DDPG agent or the optimization algorithm can generate a feasible path to the destination. Under these definitions, the optimal set is a subset of the feasible set.

\begin{figure}[htb]
\centerline{\includegraphics[height=6cm,width=9cm]{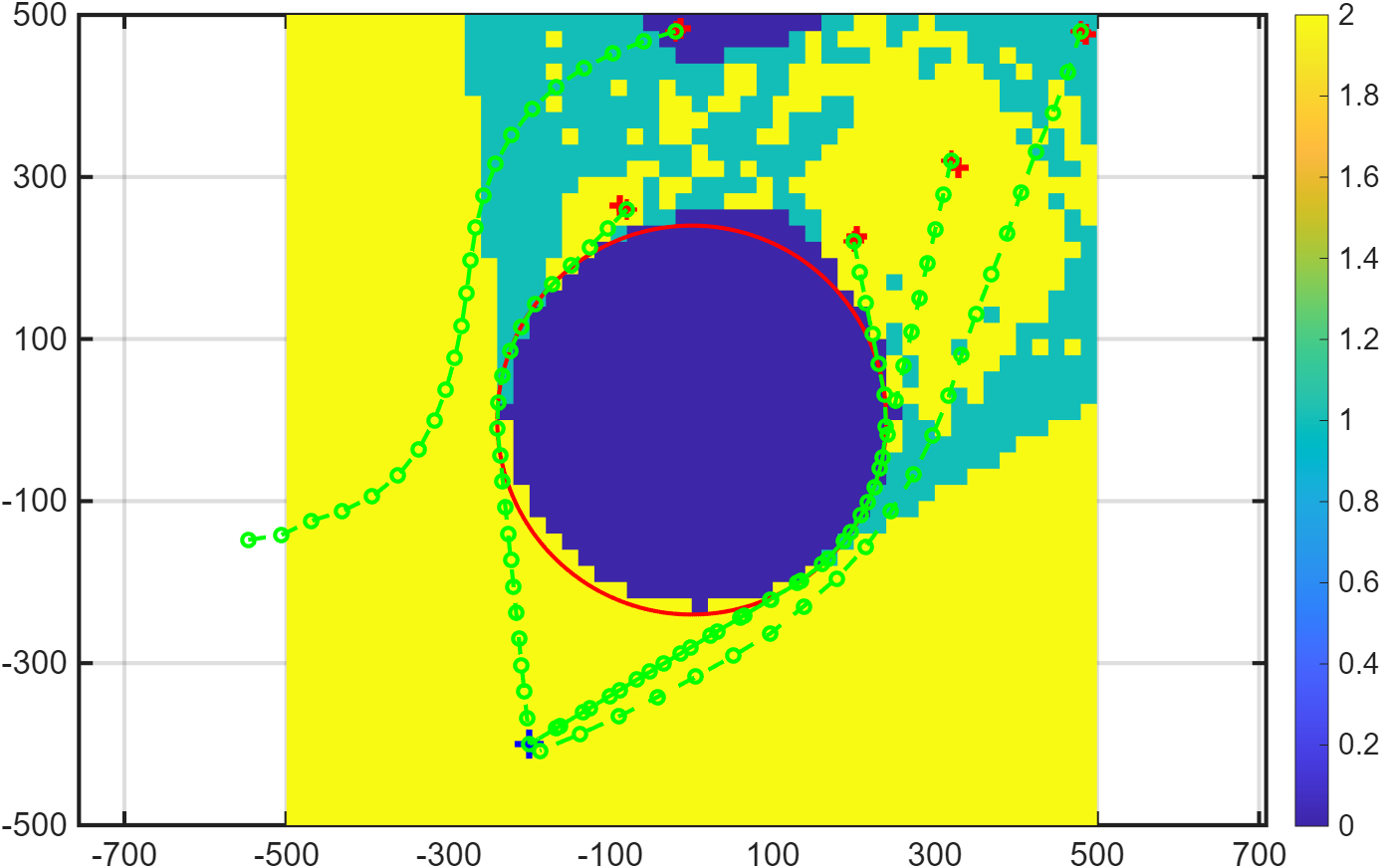}} 
\caption{Optimal and feasible sets found by DDPG-IPM implemented in Matlab fmincon.}
\label{optimalIPM}
\end{figure}


\begin{figure}[htb]
\centerline{\includegraphics[height=6cm,width=9cm]{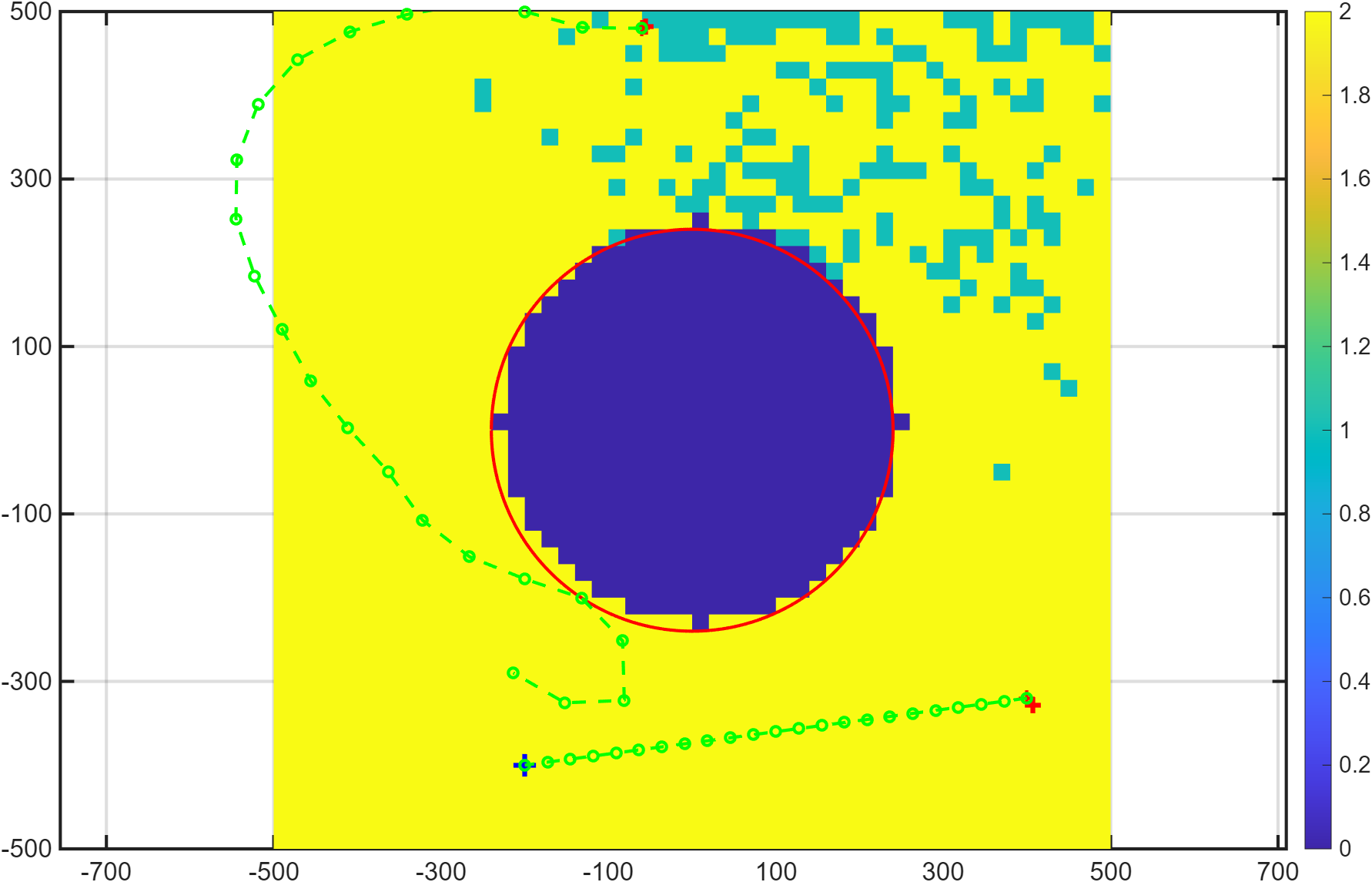}}
\caption{Optimal and feasible sets found by DDPG-SQP implemented in Matlab fmincon.}
\label{optimalSQP}
\end{figure}


In this section, we consider the path-planning problem in the presence of a single obstacle, as discussed in \cite{lyw26a}. The center of the obstacle is located at $(a,b)=(0,0)$, its radius is $R=240$, the destination is $(\bar{x}_f,\bar{y}_f)=(-200,-400)$, and the total number of line segments is $f=22$.

A path is regarded as feasible if the distance between its endpoint and the destination is less than $100$. A path is regarded as locally optimal if the corresponding segment lengths, headings, slack variables, and Lagrange multipliers satisfy the KKT conditions.

We use three algorithms to solve Problem (\ref{optP1}): the arc-search IPM algorithm (\ref{mainAlgo}), and the IPM and SQP algorithms implemented in Matlab's {\tt fmincon} function with the options {\tt interior-point} and {\tt sqp}, respectively. To ensure a fair comparison, the same initial path is used for all three algorithms for each starting point under consideration.

For a given starting point, if a feasible path obtained by DDPG is available, it is used as the initial path. Otherwise, an initial path is generated using heuristic rules. If the heuristics fail to generate an interior point (i.e., a feasible initial path satisfying the strict inequality constraints), the analytic-center method described in \cite{ck99,ye97} is employed to generate an initial path.

As a reference, the feasible set of the trained DDPG agent is shown in Figure \ref{DDPGfeasible}. In the figure, a representative starting point is marked by a red '+', the destination is marked by a blue '+', the yellow region represents the feasible set, the blue region represents the infeasible set, and a feasible path is displayed in green. The same convention is adopted in all subsequent figures.

The feasible set and the optimal set obtained by the arc-search IPM are identical, indicating that the algorithm successfully finds an optimal path for every starting point in the feasible set (100\% success rate). The result is shown in Figure \ref{arcSearch_feasible}, which demonstrates that the optimal set covers all starting points outside the obstacle. Two representative optimal paths are displayed in green.

Figure \ref{optimalIPM} shows the optimal (the area in yellow) and feasible (the area in cyan) sets obtained by the IPM implementation of Matlab's {\tt fmincon}. The blue region outside the obstacle corresponds to starting points for which {\tt fmincon} fails to compute an optimal solution. Although the Matlab IPM implementation fails to satisfy the optimality criterion for many starting points, it often terminates at solutions whose endpoints are very close to the destination. Clearly, the optimal set obtained by the Matlab IPM implementation is smaller than that obtained by the arc-search IPM. However, its feasible set is significantly larger than its optimal set.

Figure \ref{optimalSQP} shows the optimal  (the area in yellow) and feasible (the area in cyan) sets obtained by the SQP implementation of Matlab's {\tt fmincon}. Again, the optimal set is a subset of the feasible set, and both sets are smaller than those obtained by the arc-search IPM.

\subsection{Path planning in the presence of multiple obstacles}

In this section, we consider the path-planning problem in the presence of three obstacles, as discussed in \cite{lyw26a}. The centers of the obstacles are located at $(a_1,b_1)=(0,0)$, $(a_2,b_2)=(200,-400)$, and $(a_3,b_3)=(-300,200)$, with corresponding radii $R_1=240$, $R_2=150$, and $R_3=100$, respectively. The destination is $(\bar{x}_f,\bar{y}_f)=(-200,-400)$, and the total number of line segments is $f=22$. The path-planning problem is then formulated as follows:

\begin{subequations}
\begin{alignat}{2}
\min & ~~ r \label{obj3} \\
s.t.  &   ~~ \bar{x}_f=\bar{x}_0 + r \sum_{i=0}^{f-1} \cos(\theta_i),   \label{xCoor1}  \\
&   ~~  \bar{y}_f= \bar{y}_{0} +r \sum_{i=0}^{f-1} \sin(\theta_i ),   \label{yCoor1}  \\
&  ~~   \left[\bar{x}_0 + r \sum_{i=0}^{k-1} \cos(\theta_i) -a_1 \right]^2 
+ \left[\bar{y}_{0} + r \sum_{i=0}^{k-1} \sin(\theta_i) -b_1 \right]^2 -R_1^2  \geq 0, \hspace{0.1in} k=1,\ldots, f,   \label{cCons1} \\
&  ~~   \left[\bar{x}_0 + r \sum_{i=0}^{k-1} \cos(\theta_i) -a_2 \right]^2 
+ \left[\bar{y}_{0} + r \sum_{i=0}^{k-1} \sin(\theta_i) -b_2 \right]^2 -R_2^2  \geq 0, \hspace{0.1in} k=1,\ldots, f,   \label{cCons2} \\
&  ~~   \left[\bar{x}_0 + r \sum_{i=0}^{k-1} \cos(\theta_i) -a_3 \right]^2 
+ \left[\bar{y}_{0} + r \sum_{i=0}^{k-1} \sin(\theta_i) -b_3 \right]^2 -R_3^2  \geq 0, \hspace{0.1in} k=1,\ldots, f,   \label{cCons3} \\
&  ~~    \theta_k - \theta_{k-1} + 0.5 \geq 0, \hspace{0.1in} k = 1,\ldots, f-1,  \label{boundary1}  \\
&  ~~    \theta_{k-1} - \theta_{k} + 0.5 \geq 0, \hspace{0.1in} k = 1,\ldots, f-1,  \label{boundary2} \\
&  ~~  bu(k) \geq x(k) \geq bl(k),  \hspace{0.1in} k = 1,\ldots, n, \label{segLengthL1} 
\end{alignat}
\label{optP2}
\end{subequations}
where $bl(k)$ and $bu(k)$ denote the lower and upper bounds of the $k$-th variable, respectively. Specifically, $bl(1)=1$ and $bl(k)=-2\pi$ for $k\neq 1$, while $bu(1)=200$ and $bu(k)=2\pi$ for $k\neq 1$. The problem contains two equality constraints, namely (\ref{xCoor1}) and (\ref{yCoor1}); $3f$ nonlinear inequality constraints, where (\ref{cCons1}), (\ref{cCons2}), and (\ref{cCons3}) each represent $f$ nonlinear inequality constraints; $2f-2$ linear inequality constraints given by (\ref{boundary1}) and (\ref{boundary2}); and $2n$ bound constraints given by (\ref{segLengthL1}). The decision variable vector is $\x=(r,\theta_0,\theta_1,\ldots,\theta_{f-1})^T$, with $n=f+1$ variables. As in (\ref{optP}), the equality constraints can be grouped into the vector equation $\h(\x)=\mathbf{0}$, whose components are given by (\ref{xCoor1}) and (\ref{yCoor1}), while the inequality constraints can be grouped into the vector inequality $\g(\x)\geq\mathbf{0}$, whose components are given by (\ref{cCons1}), (\ref{cCons2}), (\ref{cCons3}), (\ref{boundary1}), (\ref{boundary2}), and (\ref{segLengthL1}).

To improve computational efficiency, the inequality constraints are further divided into two categories: nonlinear inequality constraints, consisting of (\ref{cCons1}), (\ref{cCons2}), and (\ref{cCons3}), and linear inequality constraints, consisting of (\ref{boundary1}), (\ref{boundary2}), and (\ref{segLengthL1}). Since the second derivatives of the linear inequality constraints are identically zero, exploiting this structure reduces the computational cost associated with Hessian evaluations. A path is regarded as locally optimal if the corresponding variables satisfy the KKT conditions for Problem (\ref{optP2}).

\begin{figure}[htb]
\centerline{\includegraphics[height=6cm,width=7cm]{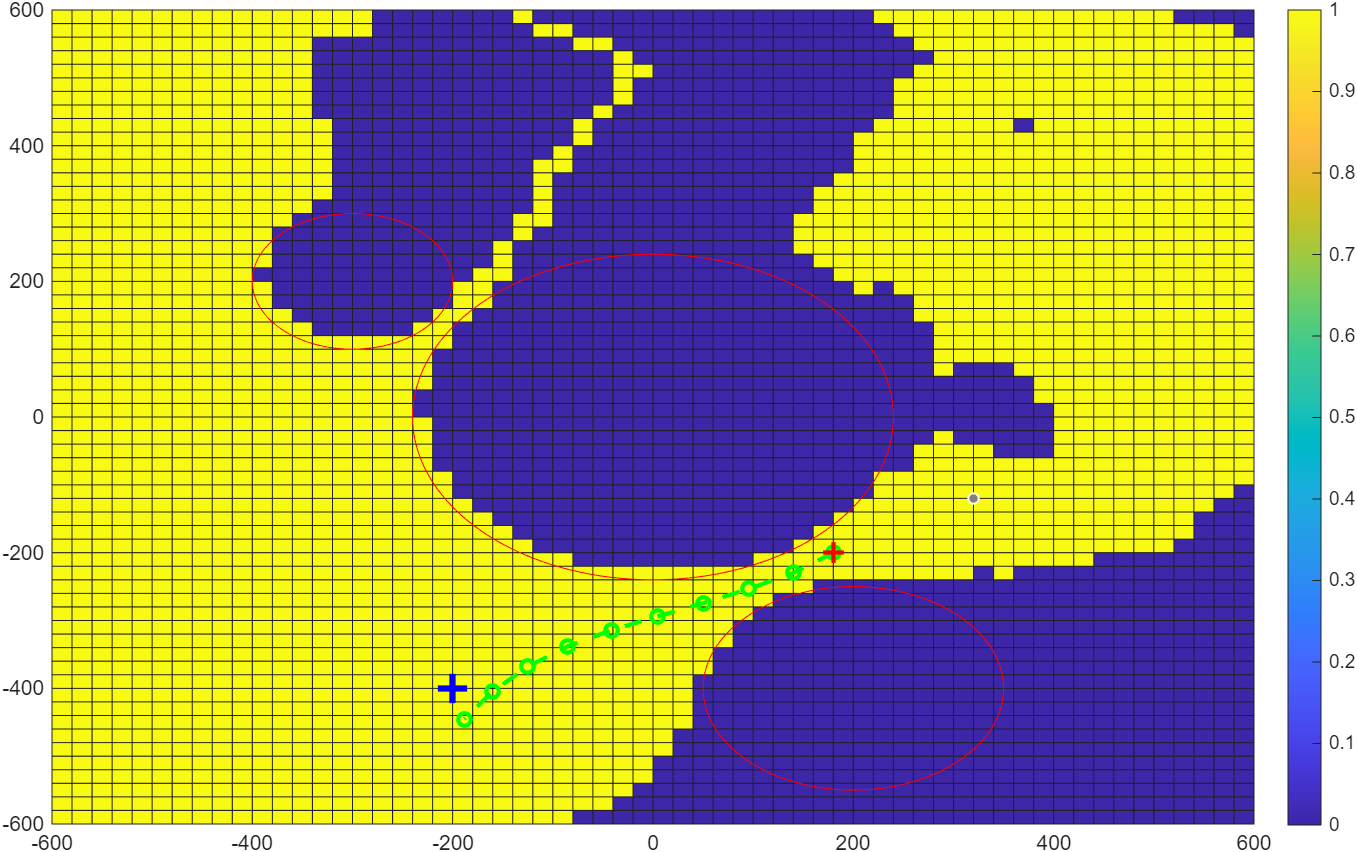}}
\caption{Feasible set found by DDPG implemented in Matlab.}
\label{threeBEZDDPG}
\end{figure}

\begin{figure}[htb]
\centerline{\includegraphics[height=6cm,width=8cm]{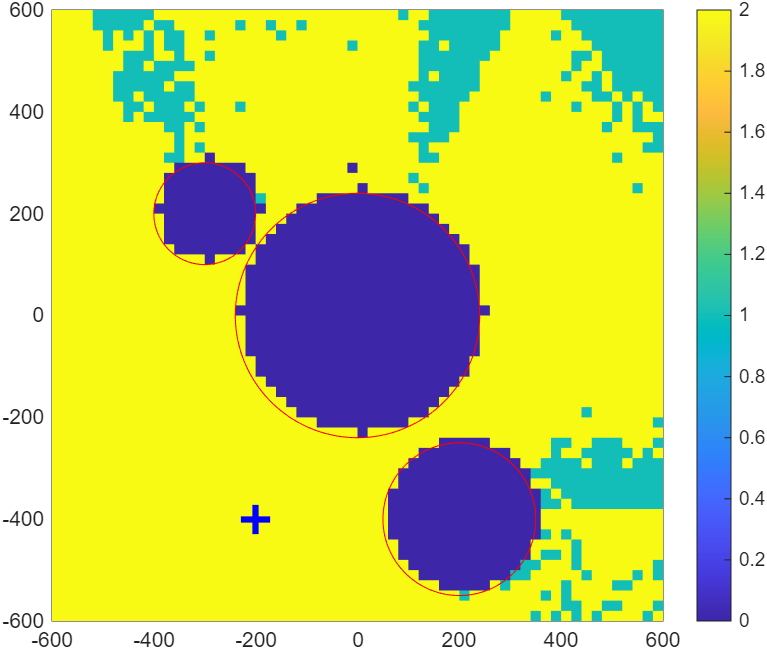}}
\caption{Feasible and optimal sets found by DDPG-arc-search implemented in Matlab.}
\label{feasible3BEZ_DDPG_arcSearch}
\end{figure}

Again, we use the reinforcement-learning algorithm DDPG to generate an initial feasible set; see \cite{lyw26a} for details. The feasible set obtained by DDPG is shown in Figure \ref{threeBEZDDPG}. Using the feasible paths generated by DDPG as initial paths, we then apply the arc-search IPM to determine both the feasible set and the optimal set. If the DDPG agent fails to generate a feasible path for a given starting point, an initial path is generated using heuristic rules. If the heuristics fail to generate an interior point, the analytic-center method described in \cite{ck99,ye97} is employed to generate an initial path. The resulting feasible set (shown in cyan) and optimal set (shown in yellow) obtained by the arc-search IPM are displayed in Figure \ref{feasible3BEZ_DDPG_arcSearch}. The figure shows that the arc-search IPM successfully finds a feasible path to the destination for almost every starting point under consideration.

\begin{figure}[htb]
\centerline{\includegraphics[height=5.5cm,width=7cm]{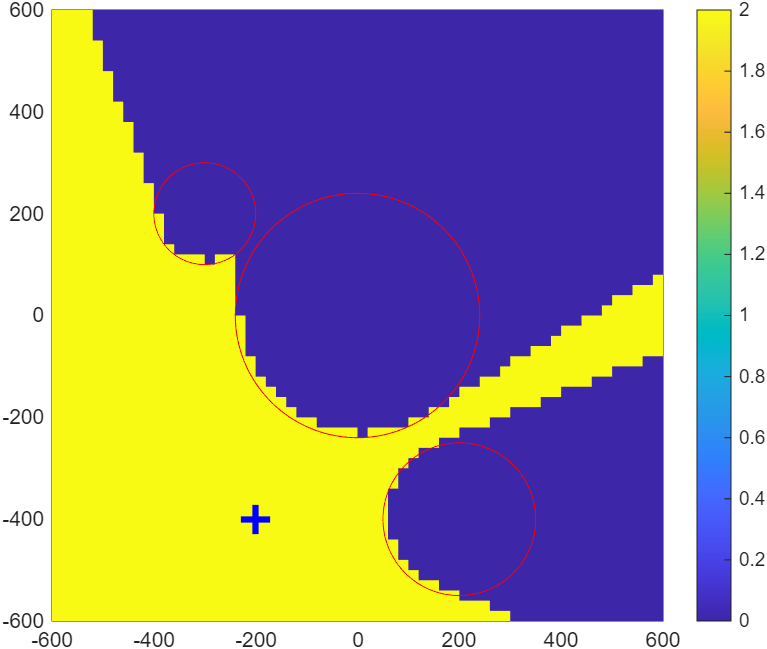}} 
\caption{Feasible and optimal sets found by DDPG-IPM implemented in Matlab.}
\label{feasible3BEZ_DDPG_IPM}
\end{figure}

\begin{figure}[htb]
\centerline{\includegraphics[height=5.5cm,width=7cm]{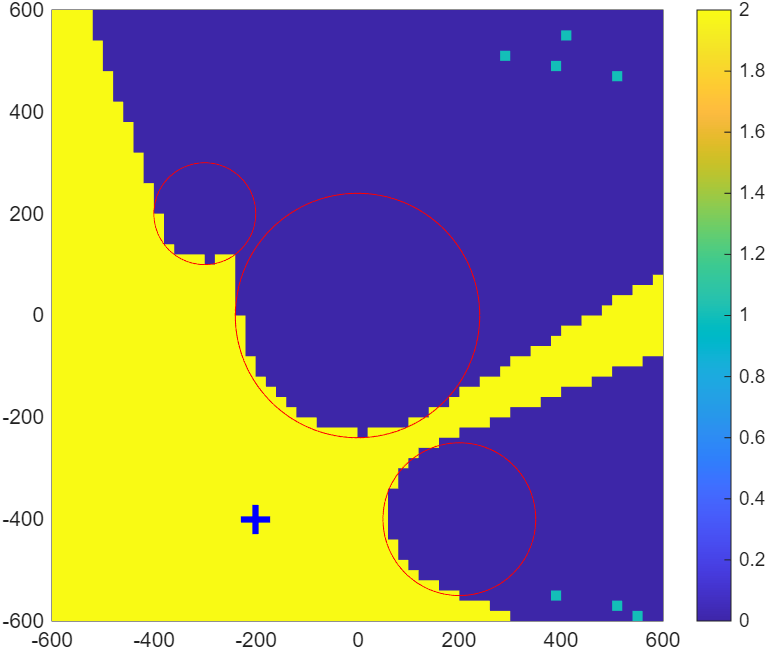}} 
\caption{Feasible and optimal sets found by DDPG-SQP implemented in Matlab.}
\label{feasible3BEZ_DDPG_SQP}
\end{figure}

The IPM and SQP implementations of Matlab’s {\tt fmincon} are also applied to this problem using the same initial paths as those used by the arc-search IPM. The corresponding feasible (shown in cyan) and optimal (shown in yellow) sets obtained by the Matlab IPM and SQP implementations are shown in Figures \ref{feasible3BEZ_DDPG_IPM} and \ref{feasible3BEZ_DDPG_SQP}, respectively. Except for starting points whose shortest path to the destination is a straight line, the Matlab IPM and SQP methods perform poorly on this problem. Even for some initially feasible paths, the Matlab IPM and SQP implementations fail to maintain feasibility while attempting to improve satisfaction of the KKT conditions.

\section{Conclusions}
\label{sec:conclusions}

In this paper, we developed a hybrid method combining deep deterministic policy gradient (DDPG) and an arc-search interior-point method (ASIPM), which leverages the strengths of each approach to mitigate the limitations of the other. Numerical simulations demonstrate that the proposed method is effective and attractive. Although the mathematical model is based on a basic engagement zone (BEZ), the proposed strategy can be readily extended to a weapon engagement zone (WEZ). Our future work will focus on extending the method to a more realistic setting with dynamically evolving engagement zones.

%
%
%
%
%
%
%
%
%
%
%

%

\end{document}